\documentclass{article}
\usepackage{graphicx}

\def\be{\begin{equation}}
\def\ee{\end{equation}}
\newcommand{\ff}[1]{{\mbox{\boldmath $#1$}}}
\def\a{\alpha}

\def\x{\ff{x}}
\begin{document}

\title{Bat Algorithm for Multi-objective Optimisation}

\author{Xin-She Yang  \\
Department of Engineering, University of Cambridge,\\
Trumpington Street,  Cambridge CB2 1PZ, UK. \\
New Address: Mathematics and Scientific Computing, \\
National Physical Laboratory, Teddington TW11 0LW, UK.}
\date{}

\maketitle

\begin{abstract}

Engineering optimization is typically multiobjective and multidisciplinary with complex constraints,
and the solution of such complex problems requires efficient optimization algorithms.
Recently, Xin-She Yang proposed a bat-inspired algorithm for solving nonlinear, global optimisation problems.
In this paper, we extend this algorithm to solve multiobjective optimisation problems.
The proposed multiobjective bat algorithm (MOBA) is first validated against a subset of test functions,
and then applied to solve multiobjective design problems such as welded beam design.
Simulation results suggest that the proposed algorithm works efficiently. \\ \\


\noindent Reference to this paper should be made as follows: \\ \\
Yang, X. S., (2011), Bat Algorithm for Multiobjective Optimization,
{\it Int. J. Bio-Inspired Computation}, Vol.~3, No.~5, pp.267-274. 

\end{abstract}


\section{Introduction}

Design optimisation in engineering often concerns multiple design objectives
under complex, highly nonlinear constraints. Different objectives often conflict
each other, and sometimes, truly optimal solutions do not exist, and some
tradeoff and approximations are often needed. Further to this complexity,
a design problem is subjected to various design constraints, limited by
design codes or standards, material properties and choice of available resources
and costs (Deb, 2001; Farina et al., 2004). Even for global optimisation problems with a single objective, if the design
functions are highly nonlinear, global optimality is not easy to reach.
Metaheuristic algorithms are very powerful in dealing with this kind of
optimization, and there are many review articles and excellent textbooks
(Coello, 1999; Deb, 2001; Isasi and Hernandez, 2004;  Yang, 2008; Talbi, 2009; Yang, 2010c).

In contrast with single objective optimization, multiobjective problems are much
difficult and complex (Coello, 1999; Floudas et al., 1999; Gong et al., 2009;
Yang and Koziel, 2010).  Firstly, no single unique solution is the best;
instead, a set of non-dominated solutions should be found in order to get a good approximation
to the true Pareto front. Secondly, even if an algorithm can find solution points
on the Pareto front, there is no guarantee that multiple Pareto points will distribute
along the front uniformly, often they do not. Thirdly, algorithms work well
for single objective optimization usually do not directly work for multiobjective problems,
unless under special circumstances such as combining multiobjectives into a single
objective using some weighted sum methods. Substantial modifications are often needed.
In addition to these difficulties, a further challenge
is how to generate solutions with enough diversity so that new solutions can
sample the search space efficiently (Talbi, 2009; Erfani and Utyuzhnikov, 2011; Yang and Koziel, 2011).

Furthermore, real-world optimization problems always involve certain degree of
uncertainty or noise. For example, materials properties for a design product
may vary significantly, an optimal design should be robust enough to allow
such inhomogeneity and also provides good choice for decision-makers or designers.
Despite these challenges, multiobjective optimization has many powerful algorithms
with many successful applications (Abbass and Sarker, 2002;
Banks et al., 2008; Deb, 2001, Farina et al., 2004; Konak et al., 2006;
Rangaiah, 2008; Marler and Arora, 2004).

In addition, metaheuristic algorithms start to emerge as
a major player for multiobjective global optimization, they often mimic the
successful characteristics in nature, especially biological systems
(Kennedy and Eberhart, 1995; Yang, 2005; Yang, 2010a; Yang, 2010b). Many new algorithms are emerging with many
important applications (Kennedy and Eberhart, 1995; Luna et al., 2007;
Osyczka and Kundu, 1995; Reyes-Sierra and Coello, 2006; Tabli, 2009;
Cui and Cai, 2009; Yang, 2010c; Zhang and Li, 2007; Yang and Deb, 2010b, Yang et al., 2011).
For example, a new cuckoo search algorithm
was developed by Xin-She Yang and Suash Deb (2009) and more detailed studies
by the same authors (Yang and Deb, 2010a) suggested that it is very efficient
for solving nonlinear engineering design problems. For a recent review of popular
metaheuristics, please refer to Yang (2011).

Recently, a new metaheuristic search algorithm, called bat algorithm (BA),
has been developed by Xin-She Yang (2010a). Preliminary studies show that
it is very promising and could outperform existing algorithms.
In this paper, we will extend BA to solve multiobjective problems and formulate
a multiobjective bat algorithm (MOBA). We will first validate it
against a subset of multiobjective test functions.
Then, we will apply it to solve design optimization problems
in engineering, such as bi-objective beam design.
Finally, we will discuss the unique features of the proposed algorithm as well as topics for further studies.

\section{Bat Behaviour and Bat Algorithm}

In order to extend the bat-inspired algorithm for single optimization to solve multiobjective
problems, let us briefly review the basics of the bat algorithm for single objective optimization.
Then, we will outline the basic ideas and steps of the proposed algorithm.

\subsection{Echolocation of Microbats}

Bats are fascinating animals. They are the only mammals with wings and
they also have advanced capability of echolocation.
It is estimated that there are about 996
different species which account for up to 20\% of all mammal
species (Altringham, 1996; Colin, 2000).
Their size ranges from the tiny bumblebee bat
(of about 1.5 to 2g) to the giant bats with wingspan of about 2 m and weight
up to about 1 kg. Microbats typically have forearm length of
about 2.2 to 11cm. Most bats uses echolocation to a certain degree;
among all the species, microbats are a famous example as microbats use
echolocation extensively while megabats do not (Richardson, 2008).

Microbats use a type of sonar, called, echolocation, to detect prey, avoid
obstacles, and locate their roosting crevices in the dark. These bats emit
a very loud sound pulse and listen for the echo that bounces back from
the surrounding objects. Their pulses vary in properties and can be
correlated with their hunting strategies, depending on the species.
Most bats use short, frequency-modulated signals to sweep through about an octave,
while others more often use constant-frequency signals for echolocation. Their signal
bandwidth varies depends on the species, and often increased by using more
harmonics.

Though each pulse only lasts a few thousandths of
a second (up to about 8 to 10 ms), however, it has a constant frequency which
is usually in the region of 25kHz to 150 kHz. The typical range of frequencies
for most bat species are in the region between 25kHz and 100kHz,
though some species can emit higher frequencies up to 150 kHz.
Each ultrasonic burst may last typically 5 to 20 ms, and microbats emit
about 10 to 20 such sound bursts every second. When hunting for prey, the rate
of pulse emission can be sped up to
about 200 pulses per  second when they fly near their prey.
Such short sound bursts imply the fantastic ability of the signal
processing power of bats. In fact, studies shows the integration
time of the bat ear is typically about 300 to 400 $\mu$s.
As the speed of sound in air is typically $v=340$ m/s, the wavelength
$\lambda$ of the ultrasonic sound bursts with a constant frequency $f$ is given by
$\lambda=v/f$,  which is in the range of 2mm to 14mm for the typical frequency range
from 25kHz to 150 kHz. Such wavelengths are in the
same order of their prey sizes.

Studies show that microbats use the time delay from the emission
and detection of the echo, the time difference between their two
ears, and the loudness variations of the echoes to build up three
dimensional scenario of the surrounding. They can detect the distance
and orientation of the target, the type of prey, and even the moving
speed of the prey such as small insects (Altringham, 1996).
Obviously, some bats have good eyesight, and most bats also have
very sensitive smell sense. In reality, they will use all the senses
as a combination to maximize the efficient detection of prey and
smooth navigation. However, here we are only interested in the echolocation
and the associated behaviour. Such echolocation behaviour of microbats
can be formulated in such a way that it can be associated with
the objective function to be optimized, and this makes it possible
to formulate new optimization algorithms.

\subsection{Bat Algorithm}

If we idealize some of the echolocation characteristics of microbats,
we can develop various bat-inspired algorithms or bat algorithms.
In the basic bat algorithm developed by Xin-She Yang (2010a), the following approximate or idealized
rules were used.

\begin{itemize}

\item[1.] All bats use echolocation to sense distance, and
they also `know' the difference between food/prey and background barriers
in some magical way;

\item[2.] Bats fly randomly with velocity $\ff{v}_i$
at position $\x_i$ with a frequency $f_{\min}$, varying wavelength $\lambda$
and loudness $A_0$ to search for prey. They can automatically adjust
the wavelength (or frequency) of their
emitted pulses and adjust the rate of pulse emission $r \in [0,1]$,
depending on the proximity of their target;

\item[3.] Although the loudness can vary in many ways, we assume that
the loudness varies from a large (positive) $A_0$ to a minimum
constant value $A_{\min}$.

\end{itemize}

Another obvious simplification is that no ray tracing is used in
estimating the time delay and three dimensional topography. Though this
might be a good feature for the application in computational geometry,
however, we will not use this feature, as it is more computationally extensive in
multidimensional cases.

In addition to these simplified assumptions, we also use the following
approximations, for simplicity. In general the frequency $f$  in
a range $[f_{\min}, f_{\max}]$ corresponds to a range
of wavelengths $[\lambda_{\min}, \lambda_{\max}]$. For example a frequency
range of [$20$kHz, $500$kHz] corresponds to a range of wavelengths
from $0.7$mm to $17$mm in reality. Obviously, we can choose the ranges freely
to suit different applications.

\subsection{Bat Motion}

For the bats in simulations, we have to define the rules how their
positions $\x_i$ and velocities $\ff{v}_i$ in a $d$-dimensional search space
are updated.
The new solutions $\x_i^{t}$ and velocities $\ff{v}_i^{t}$ at time step
$t$ are given by
\be f_i =f_{\min} + (f_{\max}-f_{\min}) \beta, \label{f-equ-150} \ee
\be \ff{v}_i^{t+1} = \ff{v}_i^{t} +  (\x_i^t - \x_*) f_i , \ee
\be \x_i^{t+1}=\x_i^{t} + \ff{v}_i^t,  \label{f-equ-250} \ee
where $\beta \in [0,1]$ is a random vector drawn from a uniform distribution.
Here $\x_*$ is the current global best
location (solution) which is located after comparing all
the solutions among all the $n$ bats at each iteration $t$.
As the product $\lambda_i f_i$ is the
velocity increment, we can use $f_i$ (or $\lambda_i$ ) to adjust the velocity change
while fixing the other factor $\lambda_i$ (or $f_i$), depending on the
type of the problem of interest.
In our implementation, we will use $f_{\min}=0$ and $f_{\max}=O(1)$, depending on the
domain size of the problem of interest. Initially, each bat is randomly assigned
a frequency which is drawn uniformly from $[f_{\min}, f_{\max}]$.

For the local search part, once a solution is selected among the current best solutions,
a new solution for each bat is generated locally using random walk
\be \x_{\rm new}=\x_{\rm old} + \ff{\epsilon}  \; A^{t}, \ee
where $\ff{\epsilon}$ is a random number vector drawn from $[-1,1]$, while $A^{t}=<\!\!A_i^{t}\!>$
is the average loudness of all the bats at this time step.

The update of the velocities
and positions of bats have some similarity to the procedure in  the standard
particle swarm optimization, as $f_i$ essentially
controls the pace and range of the movement of the swarming particles.
To a degree, BA can be considered
as a balanced combination of the standard particle swarm optimization and
the intensive local search controlled by the loudness and pulse rate.

\subsection{Loudness and Pulse Emission}

Furthermore, the loudness $A_i$ and the rate $r_i$ of pulse emission have to be
updated accordingly as the iterations proceed. As the loudness usually
decreases once a bat has found its prey, while the rate of pulse emission
increases, the loudness can be chosen as any value of convenience.
For example, we can use
$A_0=100$ and $A_{\min}=1$. For simplicity, we can also use
$A_0=1$ and $A_{\min}=0$, assuming $A_{\min}=0$ means that a bat
has just found the prey and temporarily stop emitting any sound.
Now we have
\be A_i^{t+1}=\alpha A_{i}^{t}, \;\;\;\;\; r_i^{t}= r_i^0 [1-\exp(-\gamma t)],
\label{rate-equ-50} \ee
where $\alpha$ and $\gamma$ are constants. In fact, $\alpha$ is similar
to the cooling factor of a cooling schedule in the simulated annealing (Kirkpatrick et al., 1983).
For any $0<\alpha<1$ and $\gamma>0$, we have
\be A_i^t \rightarrow 0, \;\;\; r_i^t \rightarrow r_i^0, \;\;\textrm{as} \;\;
t \rightarrow \infty. \ee
In the simplest case, we can use $\a=\gamma$,
and we have used $\a=\gamma=0.9$ in our simulations.

Preliminary studies by Yang (2010a) suggested that bat algorithm is very promising
for solving nonlinear global optimization problems. Now we extend it
to solve multiobjective optimization problems.

\section{Multiobjective Bat Algorithm}

Multiobjective optimization problems are more complicated than
single objective optimization, and we have to find and/or approximate the
optimality fronts. In addition, algorithms have to be modified to
accommodate multiobjectives properly.

\subsection{Pareto Optimality}

A solution vector $\ff{u}=(u_1, .., u_n)^T \in {\cal F}$,
is said to dominate another vector
$\ff{v}=(v_1,..., v_n)^T$ if and only if $u_i \le v_i$ for $\forall i \in \{ 1,...,n\}$
and $\exists i \in \{1,...,n \}: u_i < v_i. $
In other words, no component of $\ff{u}$ is larger than the corresponding component of $\ff{v}$,
and at least one component is smaller.
Similarly, we can define another dominance relationship $\preceq$ by
\be \ff{u} \preceq \ff{v} \Longleftrightarrow \ff{u} \prec \ff{v} \vee \ff{u}=\ff{v}. \ee
It is worth pointing out that for maximization problems, the dominance can be defined
by replacing $\prec$ with $\succ$. Therefore,
a point $\x_* \in {\cal F}$ is called a non-dominated solution if
no solution can be found that dominates it (Coello, 1999).

The Pareto front $PF$ of a multiobjective can be defined as the
set of non-dominated solutions so that
\be  PF =\{ \ff{s} \in {\cal S} \Big|
\exists \hspace{-5pt} / \; \ff{s'} \in {\cal S}: \ff{s'} \prec \ff{s} \}, \ee
or in term of the Pareto optimal set in the search space
\be  PF^* =\{ \ff{x} \in {\cal F} \Big|
\exists \hspace{-5pt} / \; \ff{x'} \in {\cal F}: \ff{f(x')} \prec \ff{f(x)} \}, \ee
where $\ff{f}=(f_1, ..., f_K)^T$.
To obtain a good approximation to Pareto front, a diverse range of solutions should be
generated using efficient techniques (Gujarathi and Babu, 2009; Konak et al., 2006).

\subsection{MOBA Algorithm}

Based on these approximations and idealization, the basic steps of the multiobjective bat algorithm
(MOBA) can be summarized as the pseudo code shown in Fig. \ref{mobafig-100}.

\begin{figure}
\begin{center}
\begin{minipage}[c]{0.9\textwidth}
\hrule \vspace{5pt}
\indent Objective functions $f_1(\x),...,f_K(\x)$,  $\x=(x_1, ..., x_d)^T$ \\
\indent Initialize  the bat population $\x_i \; (i=1,2,...,n)$ and  $\ff{v}_i$\\
\indent {\bf for} $j=1$ to $N$ (points on Pareto fronts) \\
\indent \quad Generate $K$ weights $w_k \ge 0$ so that $\sum_{k=1}^K w_k=1$ \\
\indent \quad Form a single objective $f=\sum_{k=1}^K w_k f_k$ \\
\indent \quad {\bf while} ($t<$Max number of iterations) \\
\indent \qquad Generate new solutions and update by (\ref{f-equ-150}) to (\ref{f-equ-250})  \\
\indent \qquad {\bf if } (rand $>r_i$)    \\
\indent \qquad \quad Random walk around a selected best solution \\
\indent \qquad {\bf end if} \\
\indent \qquad Generate a new solution by flying randomly \\
\indent \qquad {\bf if  } (rand $<A_i \;\; \&  \;\; f(\x_i)<f(\x_*)$) \\
\indent \qquad \quad Accept the new solutions, \\
\indent \qquad \quad and increase $r_i$ \& reduce $A_i$ \\
\indent \qquad {\bf end if} \\
\indent \qquad Rank the bats and find the current best $\x_*$ \\
\indent \quad {\bf end while} \\
\indent \quad Record $\x_*$ as a non-dominated solution \\
\indent {\bf end} \\
\indent Postprocess results and visualization
\hrule
\caption{Multiobjective bat algorithm (MOBA). \label{mobafig-100} }
\end{minipage}
\end{center}
\end{figure}

For simplicity here, we use a weighted sum to combine all objectives $f_k$
into a single objective
\be f=\sum_{k=1}^K w_k f_k, \quad \sum_{k=1}^K w_k=1. \ee
As the weights are generated randomly from a uniform distribution,
it is possible to vary the weights with sufficient diversity so that the Pareto front can be
approximated correctly.

In our simulations, we have carried out parametric studies, and we
have used $\a=\gamma=0.9$ for all simulations. The choice of parameters
requires some experimenting.
Initially, each bat should have different values of loudness and pulse emission rate,
and this can be achieved by randomization. For example, the initial loudness $A_i^{0}$
can typically be $[1, 2]$, while the initial emission rate $r_i^{0}$ can be around
zero, or any value $r_i^0 \in [0,1]$ if using (\ref{rate-equ-50}).
Their loudness and emission rates
will be updated only if the new solutions are improved, which means that these bats
are moving towards the optimal solution.

\section{Numerical Results}

\subsection{Parametric Studies}

The proposed multiobjective bat algorithm (MOBA) is implemented in Matlab, and computing time
is within a few seconds to less than a minute, depending on the problem of interest.
We have tested it using a different range of parameters such as population size ($n$),
loudness reduction $\alpha$, and pulse reduction rate $\gamma$. By varying $n=5,10,15, 20, 30, 50$
to $50, 100, 150$, $200, 250, 300, 400, 500$, $\alpha=0.5,0.6,0.7$, $0.8, 0.9, 0.95, 1$  and
$\gamma=0.5, 0.6, 0.7$, $0.8, 0.9, 0.95, 1,0$, we found that the best parameters for
most applications are: $n=25$ to $50$, $\alpha=0.7$ to $0.9$ and $\gamma=0.7$ to $0.9$.

The stopping criterion can be defined in many ways. We can either use
a given tolerance or a fixed number of iterations. From the implementation
point of view, a fixed number of iterations is not only easy to implement,
but also suitable to compare the closeness of Pareto front of different functions.
So we have set the fixed number iterations as 5000, which is sufficient
for most problems. If necessary, we can also increase it to a larger number.

In order to generate more optimal points on the Pareto front, we can
do it in two ways: increase the population size $n$ or run the program
a few more times. Through simulations, we found that to increase of $n$ typically
leads to a longer computing time than to re-run the program a few times.
This may be due to the fact that manipulations of large matrices or longer vectors
usually take longer. Another possibility is that simple restart can increase the
diversity of solutions than more intensive search for longer iterations.
So to generate 200 points using a population size $n=50$
requires to run the program 4 times, which is easily done within a few minutes.
Therefore, in all our simulations, we will use the fixed parameters: $n=50$,
$\alpha=\gamma=0.9$.

\subsection{Multiobjective Test Functions}

There are many different test functions for multobjective optimization
(Schaffer, 1985; Zhang et al., 2003; Zhang et al, 2009;
Zitzler and Thiele, 1999; Zitzler et al., 2000), but a subset of a few widely used
functions provides a wide range of diverse properties in terms Pareto front
and Pareto optimal set. To validate the proposed MOBA, we have selected a subset
of these functions with convex, non-convex and discontinuous Pareto fronts.
We also include functions with more complex Pareto sets. To be more specific
in this paper, we have tested the following four functions:

ZDT1 function with a convex front (Zitzler and Thiele 1999; Zitzler et al. 2000)
\[ f_1(x)=x_1, \quad f_2(x)=g (1-\sqrt{f_1/g}), \]
\be g=1+\frac{9 \sum_{i=2}^d x_i}{d-1}, \quad x_i \in [0,1], \; i=1,...,30, \ee
where $d$ is the number of dimensions. The Pareto-optimality is reached when $g=1$.

ZDT2 function with a non-convex front
\[ f_1(x)=x_1, \quad f_2(x) =g (1-\frac{f_1}{g})^2, \]

ZDT3 function with a discontinuous front
\[ f_1(x) =x_1, \quad f_2(x)=g \Big[1-\sqrt{\frac{f_1}{g}}-\frac{f_1}{g} \sin (10 \pi f_1) \Big], \]
where $g$ in functions ZDT2 and ZDT3 is the same as in function ZDT1. In the ZDT3 function,
$f_1$ varies from $0$ to $0.852$ and $f_2$ from $-0.773$ to $1$.

LZ4 function (Li and Zhang, 2009; Zhang and Li, 2007)
\[ f_1=x_1 +\frac{2}{|J_1|} \sum_{j \in J_1} h(u_j) \]
\be f_2=1-x_1^2 +\frac{2}{|J_2|} \sum_{j \in J_2} h(u_j), \ee
where $J_1=\{j|j$ is odd and $2 \le j \le d\}$ and $J_2 =\{ j|j$ is even and $2 \le j \le d\}$.
\[ u_j= x_j - \sin (6 \pi x_1 + \frac{j \pi}{d}), \] \[ x_1 \in [0,1], \;
x_j \in [-2,2], \; j=2,...,d, \]
and \be h(v)=\frac{|v|}{1+e^{2 |v|}}. \ee

This function has a Pareto front $f_2=1-f_1^2$ for $0 \le f_1 \le 1$ with a Pareto set
\be x_j=\sin (6 \pi x_1 + \frac{j \pi}{d}), \quad j=2,3, ..., d, \quad x_1 \in [0,1]. \ee

\begin{figure}
\centerline{\includegraphics[height=2.77in,width=4in]{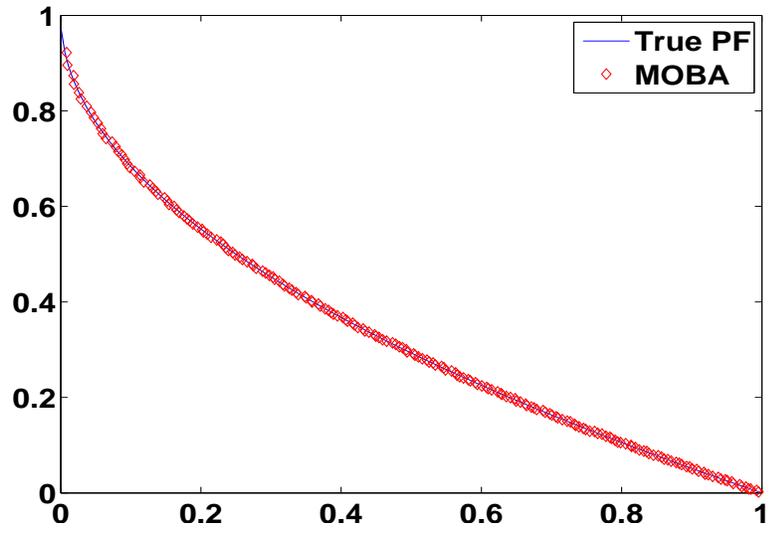}}
\caption{Estimated front and true front for ZDT1.  \label{fig-100}}
\end{figure}

After generating 200 Pareto points by MOBA, the Pareto front generated by MOBA
is compared with the true front $f_2=1-\sqrt{f_1}$ of ZDT1 (see Fig. \ref{fig-100}).
In all the rest of the figures, the vertical axis is for $f_2$ while the
horizontal axis is for $f_1$.

\begin{table}\caption{Summary of results.}
\begin{center} \begin{tabular}{|l|l|l|}
\hline
Functions & Error$_{t=2000}$ & Error$_{t=5000}$ \\
\hline
ZDT1 & 3.7E-4   & 4.5E-17  \\
ZDT2 & 2.4E-4   & 3.2E-19  \\
ZDT3 & 5.2E-5   & 1.7E-15  \\
LZ4 & 2.9E-4  & 1.2E-16 \\
\hline
\end{tabular} \end{center} \end{table}

\begin{figure}
\centerline{\includegraphics[height=2.77in,width=4in]{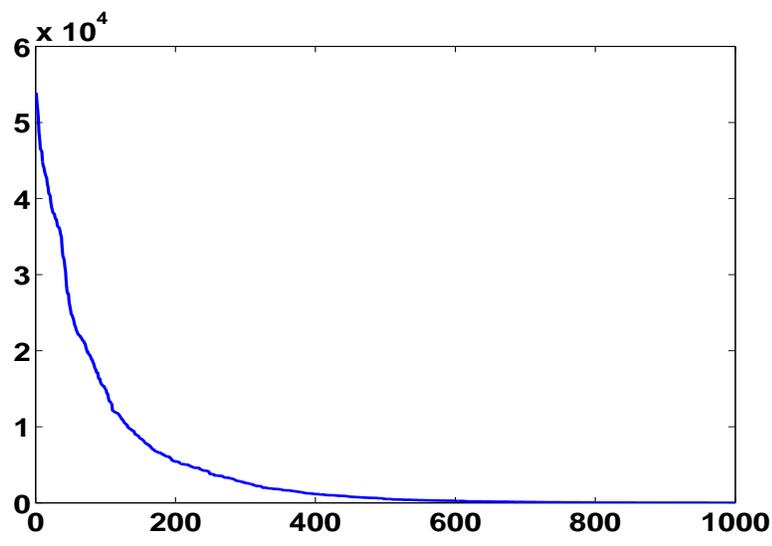}}
\caption{Convergence of the MOBA.   \label{fig-200}}
\end{figure}

\begin{figure}
\centerline{\includegraphics[height=2.77in,width=4in]{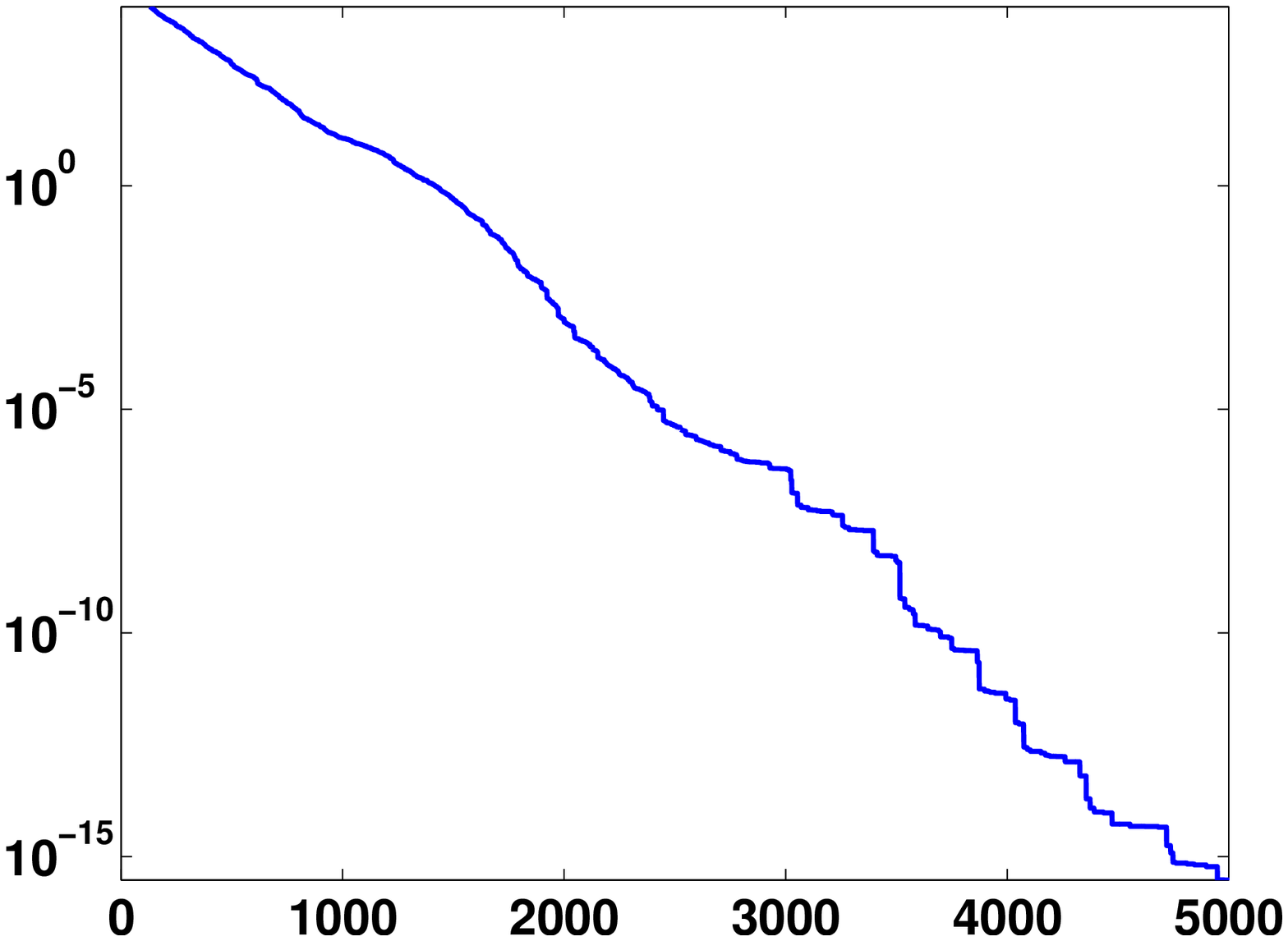} }
\caption{Exponential convergence of the MOBA.   \label{fig-300}}
\end{figure}

Let us define the distance or error between the estimate Pareto front $PF^e$ to
its correspond true front $PF^t$ as
\be E_f=||PF^e-PF^t||^2=\sum_{j=1}^N (PF_j^e-PF_t)^2, \ee
where $N$ is the number of points. The convergence property can be viewed by
following the iterations. Figs. \ref{fig-200} and \ref{fig-300} show
the exponential-like decrease of $E_f$ as the iterations proceed. The least-square distance from
the estimated front to the true front of ZDT1 for the first 1000 iterations (Fig. \ref{fig-200})
and the logarithmic scale for 5000 iterations (Fig. \ref{fig-300}).

We can see clearly that our MOBA algorithm indeed converges almost exponentially.
The results for all the functions are summarized in Table 1. We can see that
exponential convergence can be achieved in all cases.

\section{Engineering Optimization}

Design optimization, especially design of structures, has many applications
in engineering and industry. As a result, there are many different benchmarks with
detailed studies in the literature (Pham and Ghanbarzadeh, 2007; Ray and Liew, 2002;
Rangaiah, 2008).
Among the widely used benchmarks,
the welded beam design is a well-known design problem. In the rest of this paper,
we will solve this design benchmark using MOBA.

Multiobjective design of a welded beam is a classical benchmark which has been solved by
many researchers (Deb, 1999; Gong et al., 2009; Ray and Liew, 2002). The problem has four design variables: the width $w$ and length $L$
of the welded area, the depth $d$ and thickness $h$ of the main beam. The objectives are to minimize both the overall fabrication cost and the end deflection $\delta$.

The problem can be written as
\[ \textrm{minimise } \; f_1(\x)=1.10471 w^2 L + 0.04811 d h (14.0+L), \]
\be \textrm{ minimize } \; f_2=\delta, \ee
subject to
\be
\begin{array}{lll}
g_1(\x)=w -h \le 0, \vspace{2pt} \\ \vspace{3pt}
g_2(\x) =\delta(\x) - 0.25 \le 0, \\ \vspace{3pt}
g_3(\x)=\tau(\x)-13,600 \le 0, \\ \vspace{3pt}
g_4(\x)=\sigma(\x)-30,000 \le 0, \\ \vspace{3pt}
g_5(\x)=0.10471 w^2 +0.04811  h d (14+L) -5.0 \le 0, \\ \vspace{3pt}
g_6(\x)=0.125 - w \le 0, \\ \vspace{3pt}
g_7(\x)=6000 - P(\x) \le 0,
\end{array}
\ee
where
\be \begin{array}{ll}
 \sigma(\x)=\frac{504,000}{h d^2},  & Q=6000 (14+\frac{L}{2}), \\ \\
 D=\frac{1}{2} \sqrt{L^2 + (w+d)^2}, & J=\sqrt{2} \; w L [ \frac{L^2}{6} + \frac{(w+d)^2}{2}], \\ \\
 \delta=\frac{65,856}{30,000 h d^3}, &  \beta=\frac{QD}{J}, \\ \\
 \alpha=\frac{6000}{\sqrt{2} w L}, & \tau(\x)=\sqrt{\alpha^2 + \frac{\alpha \beta L}{D}+\beta^2},
 \end{array} \ee
\be P=0.61423 \times 10^6 \; \frac{d h^3}{6} (1-\frac{d \sqrt{30/48}}{28}).  \ee

The simple limits or bounds are $0.1 \le L, d \le 10$ and
$0.125 \le w, h \le 2.0$.

\begin{figure}
 \centerline{\includegraphics[height=2.5in,width=3in]{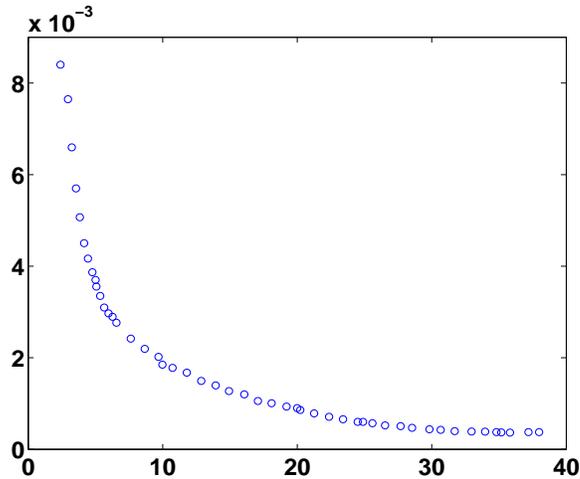} }
\caption{Pareto front for the bi-objective beam design. \label{fig-400}}
\end{figure}

By using the MOBA, we have solved this design problem. The approximate Pareto front
generated by the 50 non-dominated solutions after 1000 iterations
are shown in Fig. \ref{fig-400}. This is
consistent with the results obtained by others (Ray and Liew, 2002; Pham and Ghanbarzadeh, 2007).
In addition, the results
are more smooth with fewer iterations.

The simulations for these benchmarks and functions suggest that MOBA is a very efficient
algorithm for multiobjective optimization. It can deal with highly nonlinear problems
with complex constraints and diverse Pareto optimal sets.

\section{Conclusions}

Multiobjective optimization problems are typically very difficult to solve.
In this paper, we have successfully formulated a new algorithm for multiobjective
optimization, namely, multiobjective bat algorithm,
based on the recently developed bat algorithm. The proposed MOBA has been tested against a subset of well-chosen test functions, and
then been applied to solve design optimization benchmarks in structural engineering.
Results suggest that MOBA is an efficient multiojective optimizer.

Additional tests  and comparison of the proposed are highly needed. In the future work,
we will focus on the parametric studies for a wider range of test problems, including
discrete and mixed type of optimization problems. We will try to test the
diversity of the Pareto front it can generate so as to identify the ways to improve
this algorithm to suit a diverse range of problems.
There are a few efficient techniques to
generate diverse Pareto fronts (Erfani and Utyuzhnikov 2011), and some combination
with these techniques may improve MOBA even further.

Further research can also emphasize the performance comparison of this algorithm
with other popular methods for multiobjective optimization.
In addition, hybridization with other algorithms may also prove to be fruitful.

\end{document}